\documentstyle[11pt,amssymb]{article}
\newcommand{\beql}[1]{\begin{equation}\label{#1}}
\newcommand{\eeq}{\end{equation}}
\newcommand{\comment}[1]{}

\newcommand{\eqref}[1]{{\rm (\ref{#1})}}

\newcommand{\Abs}[1]{{\left|{#1}\right|}}

\newcommand{\Norm}[1]{{\left\|{#1}\right\|}}

\newcommand{\Qed}{\ \\\mbox{$\blacksquare$}}

\newcommand{\Set}[1]{{\left\{{#1}\right\}}}

\newcommand{\ToAppear}[1]{\raisebox{15mm}[10pt][0mm]{\makebox[0mm]%
	{\makebox[\textwidth][r]{\small #1}}}}

\newcommand{\RR}{{\Bbb R}}
\newcommand{\CC}{{\Bbb C}}
\newcommand{\ZZ}{{\Bbb Z}}
\newcommand{\NN}{{\Bbb N}}

\newcommand{\one}{{\bf 1}}

\newcommand{\inner}[2]{{\langle #1, #2 \rangle}}

\newcommand{\dens}{{\rm dens\,}}
\newcommand{\Span}{{\rm span}}

\newcommand{\supp}{{\rm supp\,}}

\newcommand{\ft}[1]{\widehat{#1}}
\newcommand{\FT}[1]{\left(#1\right)^\wedge}
\newcommand{\nozero}[1]{{#1\setminus\Set{0}}}

\newcounter{open}
\newcounter{dfn}
\def\thedfn{\arabic{dfn}}
\newenvironment{dfn}{
  \sf
  \vskip 0.10in
  \refstepcounter{dfn}
  \noindent{\bf Definition \thedfn \ }
}{\vskip 0.10in}

\newcounter{obs}
\def\theobs{\arabic{obs}}
\newenvironment{obs}{
  \sf
  \vskip 0.10in
  \refstepcounter{obs}
  \noindent{\bf Observation \theobs \ }
}{\vskip 0.10in}

\newcounter{thm}
\newenvironment{thm}{
  \sf
  \vskip 0.10in
  \refstepcounter{thm}
  \noindent{\bf Theorem\ }
}{\vskip 0.10in}

\newcounter{mysec}
\def\themysec{\arabic{mysec}}
\newcommand{\mysection}[1]{
  \vskip 0.25in
  \refstepcounter{mysec}\centerline{\large\bf \S\themysec.\ {#1}}\par
  \addcontentsline{toc}{section}{{\bf \themysec.}\ {#1}}
}

\newcounter{mysubsec}[mysec]
\def\themysubsec{\arabic{mysec}.\arabic{mysubsec}}
\newcommand{\mysubsection}[1]{
  \vskip 0.125in
  \refstepcounter{mysubsec}\noindent{\bf \themysubsec\ \ \ {#1}}\par
  \addcontentsline{toc}{subsection}{\themysubsec.\ {#1}}
}

\newtheorem{theorem}{Theorem}

\newtheorem{lemma}{Lemma}

\begin{document}

\begin{center}
{\Large \bf On the structure of multi\ToAppear{\tt PREPRINT}ple
translational tilings by\\
polygonal regions}\\
\ \\
{\sc Mihail N. Kolountzakis\footnote{Partially supported%
by the U.S. National Science Foundation,
under grant DMS 97-05775.}}\\
Department of Mathematics,\\
University of Crete,\\
Knossos Ave., 714 09 Iraklio,\\
Greece.\\
\ \\
E-mail: {\tt kolount@math.uch.gr}\\
\ \\
\small February 1999
\end{center}

\begin{abstract}
We consider polygons with the following ``pairing property'':
for each edge of the polygon there is precisely one other edge
parallel to it.
We study the problem of when such a polygon $K$ tiles the plane
multiply
when translated at the locations $\Lambda$, where $\Lambda$ is a
multiset
in the plane.
The pairing property of $K$ makes this question particularly
amenable to Fourier Analysis.
After establishing a necessary and sufficient condition for $K$ to
tile
with a given lattice $\Lambda$ (which was first found by Bolle for
the
case of convex polygons--notice that all convex polygons that tile,
necessarily have the pairing property and, therefore, our theorems
apply
to them) we move on to prove that a large class of such polygons tiles
only quasi-periodically, which for us means that $\Lambda$ must be a 
finite union of translated $2$-dimensional lattices in the plane.
For the particular case of convex polygons we show that all convex
polygons which are not parallelograms tile necessarily
quasi-periodically, if at all.
\end{abstract}

\mysection{Introduction}
\label{sec:intro}

In this paper we study multiple tilings of the plane
by translates of a polygonal region of a certain type,
the polygons with the pairing property of Definition
\ref{def:pairing} below.

\begin{dfn}\label{def:tiling} (Tiling)\\
Let $K$ be a measurable subset of $\RR^2$ of finite measure
and let $\Lambda \in \RR^2$ be
a discrete multiset (i.e., its underlying set is discrete and each
point
has finite multiplicity). We say that $K+\Lambda$ is a (translational,
multiple) \underline{tiling} of $\RR^2$, if
$$
\sum_{\lambda\in\Lambda} \one_K (x-\lambda) = w,
$$
for almost all (Lebesgue) $x\in\RR^2$, where the {\em weight}
or {\em level}
$w$ is a positive integer and $\one_K$ is the indicator function of
$K$.
\end{dfn}

\begin{dfn} \label{def:pairing}
(Polygons with the Pairing Property)\\
A polygon $K$ has the \underline{pairing property}
if for each edge $e$ there is precisely one other edge of $K$
parallel to $e$
\end{dfn}

\noindent
{\bf Remarks.}\\
1. Note that all symmetric convex polygons have the pairing
property and it is not hard to see that all convex polygons
that tile by translation are necessarily symmetric. \\
2. The polygonal regions we deal with are not assumed to be connected.

Using Fourier Analysis we study the following two problems:
(a) characterize the polygons that tile multiply with a lattice, and
(b) determine which polygons tile necessarily in a ``quasi-periodic''
manner, if they tile at all.
We restrict our attention to polygons with the pairing property.

\begin{dfn} (Quasi-periodic multisets)\\
A multiset $\Lambda \subseteq \RR^d$ is called quasi-periodic
if it is the union of finitely many $d$-dimensional lattices
(see Definition \ref{dfn:dim-lat})
in $\RR^d$.
\end{dfn}

In \S\ref{sec:fourier} we describe the
general approach to translational
tiling using the Fourier Transform of the indicator function of the
tile and in particular its zero-set.
This zero-set for polygons with the pairing property
is calculated explicitly.

In \S\ref{sec:lattice-tiling} we give a necessary and sufficient
condition (Theorem \ref{th:poly}) for a polygon $K$ with the pairing
property to
tile multiply with a lattice $\Lambda$. This has been proved
previously by Bolle for the more special case of convex polygons
(although his method might apply for the case of pairing polygons
as well) who used a combinatorial method.
Our approach is based on the calculation of
\S\ref{sec:fourier}.

In \S\ref{sec:periodic-tiling} we find a very large class
of polygons with the pairing property that tile only
in a quasi-periodic manner. In particular we show that
every convex polygon that is not a parallelogram can tile
(multiply) only in a quasi-periodic way.

\noindent
{\bf Notation.}\\
1. The Fourier Transform of a function $f \in L^1(\RR^d)$ is
normalized
as follows:
$$
\ft{f}(\xi) = \int_{\RR^d} e^{-2\pi i \inner{\xi}{x}} f(x)~dx.
$$
It is extended to tempered distributions by duality.\\
2. The action of a tempered distribution $\alpha$ on a
function $\phi$ of Schwarz class is denoted by $\alpha(\phi)$.
The Fourier Transform $\ft\alpha$ of $\alpha$ is defined by
$$
\ft\alpha(\phi) = \alpha(\ft\phi).
$$
A tempered distribution $\alpha$ is supported on a closed set $K$
if for each smooth function $\phi$ with $\supp\phi \subset K^c$ we 
have $\alpha(\phi) = 0$. The intersection of all such closed sets
$K$ is called the support of $\alpha$ and denoted by $\supp \alpha$.

\mysection{The Fourier Analytic approach.}
\label{sec:fourier}
 
\mysubsection{General}\label{sec:general}
It is easy to see that if a polygon $K$ with the pairing property
tiles multiply then for each (relevant) direction the two edges
parallel
to it necessarily have the same length. For this, suppose that $u$ is
a direction and that $e_1$ and $e_2$ are the two edges parallel to it.
Let then $\mu_u$ be the measure which is equal to arc-length on
$e_1$ and
negative arc-length on $e_2$. Suppose also that $K+\Lambda$
is a multiple tiling of
$\RR^2$. It follows then that
$$
\sum_{\lambda\in\Lambda} \mu_u ( x -\lambda)
$$
is the zero measure in $\RR^2$. This is so because each copy of edge
$e_1$ in the tiling has to be countered be some copies of edge $e_2$.
Hence the total mass of $\mu_u$ is $0$ and $e_1$ and $e_2$ have the
same length.
We can then write (here $e_1$ and $e_2$ are viewed as point-sets
in $\RR^2$ and $\tau$ as a vector)
$$
e_2 = e_1 + \tau,
$$
for some $\tau\in\RR^2$.

By the previous discussion, a polygon $K$ with the
pairing property tiles multiply with a multiset $\Lambda$ if and
only if
for each pair $e$ and $e+\tau$ of parallel edges of $K$
\beql{vanishes}
\sum_{\lambda\in\Lambda} \mu_e (x -\lambda) = 0,
\eeq
where $\mu_e$ is the measure in $\RR^2$ that is arc-length on
$e$ and negative arc-length on $e+\tau$.
Write
$$
\delta_\Lambda = \sum_{\lambda\in\Lambda} \delta_\lambda,
$$
where $\delta_a$ is a unit point mass at $a$. Thus $\delta_\Lambda$ is
locally a measure but is globally unbounded when $\Lambda$ is
infinite.
However, whenever $K+\Lambda$ is a multiple tiling, it is obvious
that $\Lambda$ cannot have more than $c R^2$ points in any disc
of radius $R$, $R>1$, ($c$ depends on $K$ and the weight of the
tiling).
This implies that $\delta_\Lambda$ is a tempered distribution and
we can take its Fourier Transform, denoted by $\ft{\delta_\Lambda}$.
Condition \eqref{vanishes} then becomes
\beql{vanishes-ft}
\ft{\mu_e}\cdot\ft{\delta_\Lambda} = 0.
\eeq

When $\Lambda$ is a lattice $\Lambda = A \ZZ^2$, where $A$ is a
$2\times2$
invertible matrix, its \underline{dual lattice} $\Lambda^*$ is defined
by
$$
\Lambda^* = \Set{x\in\RR^2:\ \inner{x}{\lambda} \in \ZZ,
 \forall\lambda\in\Lambda},
$$
and we have $\Lambda^* = A^{-\top}\ZZ^2$.
The Poisson Summation Formula then takes the form
\beql{poisson}
\ft{\delta_\Lambda} = {\det\Lambda}\cdot \delta_{\Lambda^*}.
\eeq
Since $\ft{\mu_e}$ is a continuous function
we have in this case, and whenever $\ft{\delta_\Lambda}$ is locally a
measure,
that condition \eqref{vanishes-ft} is equivalent to
\beql{supp}
\supp \ft{\delta_\Lambda} \subseteq Z(\ft{\mu_e}),
\eeq
where for every continuous function $f$ we write $Z(f)$ for the set
where it vanishes.
When $\Lambda$ is a lattice \eqref{vanishes-ft} is equivalent to
$$
\ft{\mu_e}(x) = 0, \ \ \forall x \in\Lambda^*.
$$ 

So, to check if a given polygon $K$ with the pairing property tiles
multiply $\RR^2$ with the lattice $\Lambda$,
one has to check that $\ft{\mu_e}$ vanishes
on $\Lambda^*$ for every edge $e$ of $K$.

\mysubsection{The shape of the zero-set}
Here we study the zero-set of the Fourier Transform of
the measure $\mu_e$ of \S\ref{sec:general}
and determine its structure.

We first calculate the Fourier Transform of $\mu_e$ in the particular
case when $e$ is parallel to the $x$-axis, for simplicity.
Let $\mu \in M(\RR^2)$ be the measure defined by duality by
$$
\mu(\phi) = \int_{-1/2}^{1/2} \phi(x, 0)~dx,
\ \ \forall \phi\in C(\RR^2).
$$
That is, $\mu$ is arc-length on the line segment joining the points
$(-1/2, 0)$ and $(1/2, 0)$. Calculation gives
$$
\ft{\mu}(\xi, \eta) = {\sin{\pi\xi} \over \pi\xi}.
$$
Notice that $\ft{\mu}(\xi, \eta) = 0$ is equivalent to
$\xi \in \nozero{\ZZ}$.

If $\mu_L$ is the arc-length measure on the line segment joining
$(-L/2,0)$ and $(L/2,0)$ we have
$$
\ft{\mu_L}(\xi, \eta) = {\sin{\pi L \xi} \over \pi\xi}
$$
and
$$
Z(\ft{\mu_L}) = \Set{(\xi, \eta):\ \xi\in \nozero{L^{-1}\ZZ} }.
$$
Write $\tau = (a, b)$ and let $\mu_{L,\tau}$ be the measure
which is arc-length on the segment joining $(-L/2,0)$ and $(L/2,0)$
translated by $\tau/2$ and negative arc-length on the same segment
translated by $-\tau/2$.
That is, we have
$$
\mu_{L, \tau} = \mu_L * (\delta_{\tau/2} - \delta_{-\tau/2}),
$$
and, taking Fourier Transforms, we get
$$
\ft{\mu_{L,\tau}}(\xi, \eta) = -2{\sin{\pi L \xi} \over \pi\xi}
\sin\pi(a\xi+b\eta).
$$
Define $u = {\tau\over\Abs{\tau}^2}$ and $v=(1/L, 0)$.
It follows that
$$
Z(\ft{\mu_{L,\tau}}) = (\ZZ u + \RR u^\perp) \cup (\nozero{\ZZ}v +
\RR v^\perp).
$$
This a set of straight lines of direction $u^\perp$ spaced by
$\Abs{u}$ and containing $0$ plus a similar set of lines of direction
$v^\perp$, spaced by $v$ and containing zero. However in the latter
set of parallel lines the straight line through $0$ has been removed.
We state this as a theorem for later use, formulated in a
coordinate-free way.

\begin{dfn}
(Geometric inverse of a vector)\\
The geometric inverse of a non-zero vector $u \in \RR^d$ is the vector
$$
u^* = {u \over \Abs{u}^2}.
$$
\end{dfn}
\begin{theorem}\label{th:polygon-zero-set}
Let $e$ and $e+\tau$ be two parallel line segments (translated by
$\tau$, of
magnitude and direction described by $e$, symmetric with respect to
$0$).
Let also $\mu_{e,\tau}$ be the measure which charges $e$ with its
arc-length
and $e+\tau$ with negative its arc-length.
Then
\beql{zero-shape}
Z(\ft{\mu_{e,\tau}}) = (\ZZ \tau^* + \RR \tau^{*\perp}) \cup
	(\nozero{\ZZ} e^* + \RR e^{*\perp}).
\eeq
\end{theorem}

\mysection{When does a polygon tile with a certain lattice?}
\label{sec:lattice-tiling}

The following theorem has been proved by Bolle \cite{Bo94} who used
combinatorial methods.
\begin{thm} (Bolle)\\
A convex polygon $K$, which is centrally symmetric about $0$, tiles
multiply with the lattice $\Lambda$ (for some weight $w\in\NN$)
if and only if for each edge $e$ of $K$
the following two conditions are satisfied.
\begin{itemize}
\item{(i)} In the relative interior of $e$ there is a point of
${1\over2}\Lambda$, and
\item{(ii)} If the midpoint of $e$ is not in ${1\over2}\Lambda$ then
the vector $e$ is in $\Lambda$.
\end{itemize}
\end{thm}

\noindent
{\bf Remark. } Notice that Bolle's theorem implies that all convex
polygons with vertices in $\Lambda$ tile multiply with $\Lambda$ at
some
level.

We prove the following which is easily seen to be a generalization
of Bolle's Theorem to polygons with the pairing property.
\begin{theorem}\label{th:poly}
If the polygon $K$ has the pairing property and $\Lambda$ is a
lattice in
$\RR^2$ then $K+\Lambda$ is a multiple tiling of $\RR^2$ if and only
if for each pair of edges $e$ and $e+\tau$ of $K$
\begin{itemize}
\item{(i)} $\tau \in \Lambda$, or
\item{(ii)} $e \in \Lambda$ and $\tau+\theta e \in \Lambda$, for
some $0<\theta<1$.
\end{itemize}
\end{theorem}
{\bf Proof of Theorem \ref{th:poly}.}
Once again we simplify matters and take the edge $e$ to be parallel
to the $x$-axis and follow the notation of \S\ref{sec:general}.

For an arbitrary non-zero vector $w\in\RR^2$ define the group
$$
G(w) = \ZZ w + \RR w^\perp,
$$
which is a set of straight lines in $\RR^2$ of direction $w^\perp$
spaced regularly at distance $\Abs{w}$.
It follows that 
$$
Z(\ft{\mu_{L,\tau}}) \subseteq G(u) \cup G(v).
$$
 From Theorem \ref{th:polygon-zero-set} it follows that
$\Lambda^* \subseteq Z(\ft{\mu_{L,\tau}})$ which implies that
$\Lambda^* \subset G(u)$ or $\Lambda^* \subset G(v)$.

This is a consequence of the following.
\begin{obs}\label{obs:groups}
If $G, H, K$ are groups and $G\subseteq H\cup K$
then $G\subseteq H$ or $G\subseteq K$.
\end{obs}
For, if $a \in G\setminus K$ and $b\in G\setminus H$, then $a\cdot
b \in H$, say, which implies $b\in H$, a contradiction.

So we have the two alternatives
\begin{enumerate}
\item\label{alt:1}
$\Lambda^* \subset G(u)$,
\item\label{alt:2}
$\Lambda^* \subset G(v)$.
\end{enumerate}
However, since not all of $G(v)$ is in $Z(\ft{\mu_{L,\tau}})$, if
alternative
\ref{alt:2} holds and alternative \ref{alt:1} does not, it follows
that
\beql{tmp1}
\Lambda^* \subseteq \Span_\ZZ\Set{v, w},
\eeq
where $w$ is the smallest (in length) multiple of $v^\perp$ which is
in $G(u)$,
i.e.,
$$
w = (0, 1/b).
$$
We have that \eqref{tmp1} is equivalent to
$$
\Lambda \supseteq \left( \Span_\ZZ\Set{v, w} \right)^* 
 = \ZZ (L, 0) + \ZZ (0, b),
$$
which is in turn equivalent to
$$
(L,0) \in \Lambda\ \ \ \mbox{and}\ \ \ (0,b) \in \Lambda.
$$
Notice also that
$$
\Lambda^* \subseteq G(u) \Longleftrightarrow
\Lambda \supseteq {G(u)}^* \Longleftrightarrow
\Lambda \ni {u\over \Abs{u}^2} = \tau.
$$

We have therefore proved the following lemma.
\begin{lemma}\label{lm:lm1}
If $\Lambda$ is a lattice, $u = {(a,b)\over a^2+b^2}$ and
$v = (L, 0)$, then
$$
\Lambda^* \subset (\ZZ u + \RR u^\perp) \cup (\nozero{\ZZ}v + \RR
v^\perp)
$$
if and only if
\begin{enumerate}
\item
$(a, b) \in \Lambda$, or
\item
$(L, 0) \in \Lambda$ and $(0,b) \in \Lambda$.
\end{enumerate}
\end{lemma}

Allowing for a general linear transformation, let $\tau, e \in
\RR^2$, and let
$\mu_{e,\tau}$ be the measure that ``charges''
with its arc-length the line segment $e$
translated so that its midpoint is at $\tau/2$ and charges with
negative its
arc-length the line segment $e$ with its midpoint at $-\tau/2$.
We have proved the following:
\beql{general-position}
\Lambda^* \subset Z(\ft{\mu_{e,\tau}}) \Longleftrightarrow
\left\{\begin{array}{l}
\tau \in \Lambda,\ \ \mbox{or}\\
e \in \Lambda\ \ \mbox{and}\ \ \tau+\theta e \in \Lambda,
 \ \ \mbox{for some $0<\theta<1$}.
\end{array}\right.
\eeq
This completes the proof of Theorem \ref{th:poly}.
\Qed

\mysection{Polygons that tile only quasi-periodically}
\label{sec:periodic-tiling}

\mysubsection{Meyer's theorem}\label{sec:meyer}
We now deal with the following question: which polygons with
the pairing property admit only quasi-periodic multiple tilings.
The main tool here, as was in \cite{KL96}, is the idempotent theorem
of P.J. Cohen
for general locally compact abelian groups, in the form of the
following
theorem of Y. Meyer \cite{M70}.

\begin{dfn}\label{dfn:coset-ring}
(The coset ring)\\
The \underline{coset ring} of an abelian group $G$ is the smallest
collection of subsets of $G$ which is closed under finite unions,
finite intersections and complements (that is, the smallest
\underline{ring} of subsets of $G$) and which contains
all cosets of $G$
\end{dfn}

\noindent
{\bf Remark.} When the group is equipped with a topology one
usually only demands that the open cosets of $G$
are in the coset ring, but we take all cosets in our definition.
\begin{thm}
(Meyer)\\
Let $\Lambda \subseteq \RR^d$ be a discrete set and $\delta_\Lambda$
be the Radon measure
$$
\delta_\Lambda = \sum_{\lambda \in \Lambda} c_\lambda \delta_\lambda,
\ \ c_\lambda\in S,
$$
where $S\subseteq\nozero{\CC}$ is a \underline{finite} set.
Suppose that $\delta_\Lambda$ is tempered,
and that $\ft{\delta_\Lambda}$ is a Radon measure on $\RR^d$ which
satisfies
\beql{growth}
\Abs{\ft{\delta_\Lambda}} ([-R,R]^d) \le C R^d, ~ \mbox{as} ~ R
\rightarrow
  \infty,
\eeq
where $C>0$ is a constant.
Then, for each $s\in S$,
the set
$$
\Lambda_s = \Set{\lambda\in\Lambda:\ c_\lambda = s}
$$
is in the coset ring of $\RR^d$.
\end{thm}

A proof of Meyer's theorem for $d=1$ can be found in \cite{KL96}.
The proof works verbatim for all $d$.

\mysubsection{Discrete elements of the coset ring}
\label{sec:coset-ring}

In this section we determine the structure of the discrete elements
of the coset ring of $\RR^d$.

In dimension $d=1$ we have the following characterization of the
discrete elements of the coset ring of $\RR$, due to Rosenthal
\cite{R66}.
\begin{thm} (Rosenthal)\\
The elements of the coset ring of $\RR$ which are discrete in the
usual topology of $\RR$ are precisely the sets of the form
\beql{ap-unions}
F \bigtriangleup \bigcup_{j=1}^J ( \alpha_j \ZZ + \beta_j ) ~,
\eeq
where $F \subseteq \RR$ is finite, $\alpha_j > 0$ and $\beta_j \in
\RR$
($\bigtriangleup$ denotes symmetric difference).
\end{thm}
Rosenthal's proof does not extend to dimension $d\ge 2$.
Since we need to know what kind of sets the elements of the
coset ring of $\RR^2$ are, we prove the following general
theorem.
\begin{theorem}\label{th:coset-ring}
Let $G$ be a topological abelian group and let ${\cal R}$ be
the least ring of sets which contains the discrete cosets of
$G$.
Then ${\cal R}$ contains all discrete elements of the coset ring of
$G$.
\end{theorem}
In other words, a discrete element of the coset ring can always
be written as a finite union of sets of the type
\beql{intersections1}
A_1 \cap \cdots \cap A_m \cap B_1^c \cap \cdots \cap B_n^c,
\eeq
where the $A_i$ and $B_i$ are \underline{discrete} cosets of $G$.
And, observing that the intersection of any two cosets is a coset,
we may rewrite \eqref{intersections1} as
\beql{intersections2}
A \cap B_1^c \cap \cdots \cap B_n^c,
\eeq
where $A$ and all $B_i$ are discrete cosets.

We need the following lemma.
\begin{lemma}\label{lm:group}
Suppose that $A$ is a non-discrete topological abelian group,
$F \subset A$ is discrete and $B_1,\ldots,B_m$ are cosets
in $A$ disjoint from $F$.
Then
\beql{union-of-cosets}
A = F \cup B_1 \cup \cdots \cup B_n
\eeq
implies that $F=\emptyset$.
This remains true if $A$ is a coset in a larger group.
\end{lemma}
{\bf Proof of Lemma \ref{lm:group}.}
Write $B_i = x_i + G_i$ and let $k$ be the number of different
subgroups
$G_i$ appearing in \eqref{union-of-cosets}.
We do induction on $k$.
Notice that the group $G_1$ may be assumed to be non-discrete, by the
non-discreteness of $A$.

When $k=1$ the theorem is true as then $F$ is a union of cosets of
$G_1$ and cannot be discrete unless it is empty.
(Here is where the disjointness of $F$ from the $B_i$ is used.)

Assume the theorem true for $k\le n$ and suppose that precisely $n+1$
groups appear in \eqref{union-of-cosets} and that $F \neq \emptyset$.
Assume that the $G_1$-cosets in \eqref{union-of-cosets} are
$$
x_1+G_1,\ldots,x_r+G_1,
$$
and let $y \in F$.
We then have
$$
y+G_1 \subseteq F \cup (X_2 + G_2) \cup \cdots \cup (X_{n+1}+G_{n+1}),
$$
with all sets $X_i$, $i=2,\ldots,n+1$, being finite.
Hence
\begin{eqnarray*}
G_1 &\subseteq& (-y +F) \cup (-y+X_2+G_2) \cup \cdots \cup (-
y+X_{n+1}+G_{n+1}) \\
 &=& F' \cup (X_2' + G_2) \cup \cdots \cup (X_{n+1}' + G_{n+1}),
\end{eqnarray*}
with $F' = -y+F$, $X_i' = -y + X_i$.

Furthermore, one may take $X_i' \subset G_1$, $i=2,\ldots,n+1$
(possibly empty), to get
$$
G_1 \subseteq (F' \cap G_1) \cup (X_2' + G_2\cap G_1) \cup \cdots
 \cup (X_{n+1}' + G_{n+1}\cap G_1).
$$
Since $y\in F$ we have that $F'\cap G_1 \ni 0$ (hence it is non-empty)
and
$$
(F'\cap G_1) \cap (X_i'+G_i\cap G_1) = \emptyset,\ \ \ i=2,\ldots,n+1.
$$
By the induction hypothesis we get a contradiction.
\Qed

\noindent
{\bf Proof of Theorem \ref{th:coset-ring}.}
By Lemma \ref{lm:group}, if $A$ is non-discrete then
$A \cap B_1^c \cap \cdots \cap B_n^c$ is either non-discrete or empty.
Hence a finite union of such sets can only be discrete
if all participating $A$'s are discrete.
Rewrite then
$$
A \cap B_1^c \cap \cdots \cap B_n^c =
 A \cap (B_1 \cap A)^c \cap \cdots \cap (B_n \cap A)^c
$$
so as to have the arbitrary discrete element of the coset ring
made up with finitely many operations from discrete cosets.
\Qed

\begin{dfn} \label{dfn:dim-lat} (Dimension, lattices)\\
The \underline{dimension} of a set $A \subseteq \RR^d$ is the
dimension
of the smallest translated subspace of $\RR^d$ that contains $A$.
A \underline{lattice} is a discrete subgroup of $\RR^d$.
\end{dfn}

\noindent
{\bf Remark.} It is well known that all $k$-dimensional lattices
in $\RR^d$ are of the form $A \ZZ^k$, where $A$ is a $d\times k$
real matrix of rank $k$.

\begin{theorem}\label{th:discrete-lattice}
Let $C = A \cap B_1^c \cap \cdots \cap B_n^c$, with $A, B_i$
being discrete cosets of $\RR^d$.
Then $C$ may be written as a finite (possibly empty) union of sets
of the type
\beql{block}
K \cap L_1^c \cap \cdots \cap L_m^c,\ \ \ L_i \subseteq K \subseteq A,
 \ \ \ m\ge 0,
\eeq
where the $K, L_i$ are discrete cosets and, when $C$ is not empty,
$$
\dim L_i < \dim K = \dim A = \dim C.
$$
\end{theorem}

\begin{obs}\label{obs:overlaps}
If $A$ and $B$ are discrete cosets in $\RR^d$ with
$\dim A = \dim B = \dim A\cap B$ then $A\cap B^c$ is a finite
(possibly empty) union of disjoint cosets of $A\cap B$ and, therefore,
$\dim A\cap B^c = \dim A$, except when $A \cap B^c = \emptyset$.
Hence $A$ and $B$ can each be written as a finite disjoint union of
translates of $A \cap B$.
\end{obs}

\noindent
{\bf Proof of Theorem \ref{th:discrete-lattice}:}
Notice that
$$
C = A \cap (B_1\cap A)^c \cap \cdots \cap (B_n \cap A)^c.
$$
Let
$$
\alpha = \dim A = \dim B_1 \cap A = \cdots = \dim B_r \cap A
$$
and
$\dim B_i \cap A < \alpha$ for $i > r \ge 0$.
Let
$$
C' = A \cap (B_1 \cap A)^c \cap \cdots \cap (B_r\cap A)^c.
$$
By induction on $r\ge 0$ we prove that $C'$ is a finite union of
sets of
type
\eqref{block}.
For $r=0$ this is obvious.
If it is true for $r-1$ then $C'$ is a finite union of sets of type
$$
K \cap L_1^c \cap \cdots \cap L_m^c \cap (B_r \cap A)^c,
$$
with $\alpha = \dim K > \dim L_i$, $i=1,\ldots,m$.
Each of these sets falls into one of two categories:

\noindent
\underline{Category 1:} $\dim K \cap (B_r \cap A) = \alpha$.\\
Then, by Observation \ref{obs:overlaps}
above, $K \cap (B_r \cap A)^c$ is a finite
union of cosets $K_1,\ldots,K_s$ of dimension $\alpha$ and hence
$C'$ is a finite union of $K_i \cap L_1^c \cap \cdots \cap L_m^c$.

\noindent
\underline{Category 2:} $\dim K \cap (B_r \cap A) < \alpha$.\\
Then
$$
K \cap L_1^c \cap \cdots \cap L_m^c \cap (B_r \cap A)^c
$$
is already of the desired form.
\Qed

From Theorems \ref{th:coset-ring} and \ref{th:discrete-lattice}
it follows for $d=2$ that every discrete element $S$ of the coset
ring of $\RR^2$ may be written as 
\beql{form}
S ~ = ~ \left( \bigcup_{j=1}^J A_j ~ \setminus ~ (B_1^{(j)} \cup 
\cdots \cup
 B_{n_j}^{(j)})\right)~ \cup ~ \bigcup_{l=1}^L L_l ~ \bigtriangleup
~ F,
\eeq
where $A_1,\ldots,A_J$ are $2$-dimensional translated lattices,
$L_l$ and
$B_i^{(j)}$ are
$1$-dimensional translated lattices and $F$ is a finite set ($J, L
\ge 0$).
And, repeatedly using Observation \ref{obs:overlaps}, the lattices
$A_j$ may be assumed to be have pairwise intersections of
dimension at most $1$.

\mysubsection{Purely discrete Fourier Transform}
\label{sec:purely-discrete}

\begin{dfn}
(Uniform density)\\
A multiset $\Lambda\subseteq\RR^d$ has asymptotic density
$\rho$ if
$$
\lim_{R\to\infty} {\Abs{\Lambda \cap B_R(x)} \over \Abs{B_R(x)}} \to
\rho
$$
uniformly in $x\in\RR^d$.\\
We say that $\Lambda$ has (uniformly) bounded density if the fraction
above is bounded by a constant $\rho$ uniformly for $x\in\RR$ and
$R>1$.
We say then that $\Lambda$ has density uniformly bounded by $\rho$.
\end{dfn}

Assume that $\Lambda \subset \RR^2$ is a discrete multiset
of bounded density which
satisfies the assumptions of Meyer's Theorem
(if we write $c_\lambda$ for the multiplicity of $\lambda\in\Lambda$).
Then, if $\Lambda_k$ is the subset of $\Lambda$ of multiplicity $k$,
$\Lambda_k$ is a discrete element of
the coset ring and is of the form \eqref{form}.

Assume now in addition that $\ft{\delta_\Lambda}$ has discrete
support.
We shall prove that all sets $F$, $L_l$ and $B_i^{(j)}$ are empty in
\eqref{form} and so
$$
\Lambda = \bigcup_{j=1}^J A_j,
$$
where the $A_i$ are translated $2$-dimensional lattices in $\RR^2$.

One can easily show that whenever $\Omega\subseteq\RR^d$ of finite
measure
tiles with $\Lambda$ at level $w$ then $\Lambda$ has density
$w/\Abs{\Omega}$.

\begin{theorem}\label{th:at-zero}
Suppose that $\Lambda \in \RR^d$ is a multiset with density $\rho$,
$\delta_\Lambda = \sum_{\lambda\in\Lambda} \delta_\lambda$,
and that $\ft{\delta_\Lambda}$ is a measure in a neighborhood of
$0$.
Then $\ft{\delta_\Lambda}(\Set{0}) = \rho$.
\end{theorem}
{\bf Proof of Theorem \ref{th:at-zero}.}
Take $\phi \in C^\infty$ of compact support with $\phi(0) = 1$.
We have
\begin{eqnarray*}
\ft{\delta_\Lambda}(\Set{0}) &=&
 \lim_{t\to\infty} \ft{\delta_\Lambda}(\phi(t x))\\
&=& \lim_{t\to\infty} \delta_\Lambda(t^{-d}\ft\phi(\xi/t))\\
&=& \lim_{t\to\infty} t^{-d}
\sum_{\lambda\in\Lambda}\ft\phi(\lambda/t)\\
&=& \lim_{t\to\infty} \sum_{n\in\ZZ^d} \sum_{\lambda\in Q_n} t^{-d}
\ft\phi(\lambda/t)
\end{eqnarray*}
where, for fixed and large $T>0$,
$$
Q_n = [0,T)^d + T n,\ \ \ \ n \in \ZZ^d.
$$
Since $\Lambda$ has density $\rho$ it follows that for each
$\epsilon>0$
we can choose $T$ large enough so that for all $n$
$$
\Abs{\Lambda \cap Q_n} = \rho \Abs{Q_n}(1 + \delta_n),
$$
with $\Abs{\delta_n} \le \epsilon$.
For each $n$ and $\lambda\in Q_n$ we have
$$
\ft\phi(\lambda / t) = \ft\phi(T n / t) + r_\lambda
$$
with $\Abs{r_\lambda} \le C T t^{-1} 
 \Norm{\nabla\ft\phi}_{L^\infty(t^{-1}Q_n)}$.
Hence
\begin{eqnarray*}
\ft{\delta_\Lambda}(\Set{0}) &=&
 \lim_{t\to\infty} \sum_{n\in\ZZ^d} t^{-d} \sum_{\lambda\in Q_n}
   (\ft{\phi}(Tn/t) + r_\lambda) \\
&=& \lim_{t\to\infty} \sum_{n\in\ZZ^d} t^{-d} 
	\rho\Abs{Q_n}(1+\delta_n)\ft\phi(Tn/t) + \\
&\ & \lim_{t\to\infty} \sum_{n\in\ZZ^d} 
		t^{-d} \sum_{\lambda\in Q_n} r_\lambda \\
&=& \lim_{t\to\infty}S_1 + \lim_{t\to\infty}S_2.
\end{eqnarray*}
We have
\beql{tmp2}
\Abs{S_1 - \sum_n t^{-d} \rho \Abs{Q_n} \ft\phi(Tn/t)} \le
  \epsilon \sum_n t^{-d} \rho \Abs{Q_n} \Abs{\ft\phi(Tn/t)}
\eeq
The first sum in \eqref{tmp2} is a Riemann sum for 
$\rho\int_{\RR^d}\ft\phi = \rho$ and the second is
a Riemann sum for $\rho\int_{\RR^d}\Abs{\ft\phi}<\infty$.

For $S_2$ we have
\begin{eqnarray*}
\Abs{S_2} &\le&
 C \sum_{n\in\ZZ^d} t^{-d} \rho \Abs{Q_n} (1+\delta_n) T t^{-1}
   \Norm{\nabla\ft\phi}_{L^\infty(t^{-1}Q_n)} \\
&\le& \rho C T t^{-1} \sum_{n\in\ZZ^d} t^{-d}\Abs{Q_n}
   \Norm{\nabla\ft\phi}_{L^\infty(t^{-1}Q_n)}.
\end{eqnarray*}
The sum above is a Riemann sum for $\int_{\RR^d} \Abs{\nabla\ft\phi}$,
which is finite, hence $\lim_{t\to\infty} S_2 = 0$.

Since $\epsilon$ is arbitrary the proof is complete.
\Qed

\noindent
{\bf Remark:}
The same proof as that of Theorem \ref{th:at-zero} shows that,
if
$$
\mu = \sum_{\lambda\in\Lambda} c_\lambda \delta_\lambda,
$$
with $\Abs{c_\lambda} \le C$, $\Lambda$ is of density $0$ and
the tempered distribution $\ft\mu$ is locally a measure
in the neighborhood of some point $a\in\RR^2$, then we 
have $\ft\mu(\Set{a}) = 0$.

\begin{theorem}\label{th:periodic-set}
Suppose that $\Lambda \subset \RR^2$ is a uniformly discrete
multiset and that
$$
\ft{\delta_\Lambda} = \FT{\sum_{\lambda\in\Lambda} \delta_\lambda}
$$
is locally a measure with
$$
\Abs{\ft{\delta_\Lambda}}(B_R(0)) \le C R^2,
$$
for some positive constant $C$.
Assume also that $\ft{\delta_\Lambda}$ has discrete support.
Then $\Lambda$ is a finite union of translated lattices.
\end{theorem}
{\bf Proof of Theorem \ref{th:periodic-set}.}
Define the sets (not multisets)
$$
\Lambda_k = \Set{\lambda \in \Lambda:\ \mbox{$\lambda$
 has multiplicity $k$}}.
$$
By Meyer's Theorem (applied for the base set of the multiset $\Lambda$
with the coefficients $c_\lambda$ equal to the corresponding
multiplicities) each of the $\Lambda_k$ is in the coset ring
of $\RR^2$ and, being discrete, is of the type \eqref{form}.

We may thus write
\beql{identity}
\Lambda_k = A \bigtriangleup B,
\eeq
with $A = \bigcup_{j=1}^J A_j$, where the $2$-dimensional
translated lattices $A_j$ have pairwise intersections of
dimension at most $1$, and $\dens B = 0$.

Hence
$$
\delta_{\Lambda_k} = \sum_{j=1}^J \delta_{A_j} + \mu,
$$
where $\mu = \sum_{f \in F} c_f \delta_f$, $\dens F = 0$ and
$\Abs{c_f} \le C(J)$.
The set $F$ consists of $B$ and all points contained in at least
two of the $A_j$.

Combining for all $k$, and reusing the symbols
$A_j$, $\mu$ and $F$ we get
$$
\delta_{\Lambda} = \sum_{j=1}^J \delta_{A_j} + \mu.
$$
But $\ft{\delta_\Lambda}$ and $\sum_{j=1}^J \ft{\delta_{A_j}}$ are
both (by the assumption and the Poisson Summation Formula) discrete
measures, and so is therefore $\ft\mu$.
However $\dens F = 0$ and the boundedness of the coefficients $c_f$
implies that $\ft\mu$ has no point masses
(see the Remark after the proof of Theorem \ref{th:at-zero}),
which means that $\ft\mu = 0$
and so is $\mu$.
Hence $\delta_{\Lambda} = \sum_{j=1}^J \delta_{A_j}$, or
$$
\Lambda = \bigcup_{j=1}^J A_j,\ \ \ \mbox{as multisets}.
$$
\Qed

Finally, we show that discrete support for
$\ft{\delta_\Lambda}$ implies that $\ft{\delta_\Lambda}$
is locally a measure.
\begin{theorem}\label{th:point-implies-measure}
Suppose that the multiset $\Lambda\subset\RR^d$ has density
uniformly bounded by $\rho$
and that, for some point $a\in\RR^d$ and $R>0$,
$$
\supp{\ft{\delta_\Lambda}} \cap B_R(a) = \Set{a}.
$$
Then, in $B_R(a)$, we have $\ft{\delta_\Lambda} = w \delta_a$, for
some $w \in \CC$ with $\Abs{w} \le \rho$.
\end{theorem}  
{\bf Proof of Theorem \ref{th:point-implies-measure}.}
It is well known that the only tempered distributions supported
at a point $a$ are finite linear combinations of the derivatives
of $\delta_a$.
So we may assume that, for $\phi\in C^\infty(B_R(a))$,
\beql{point-distr}
\ft{\delta_\Lambda}(\phi) =
 \sum_\alpha c_\alpha (D^\alpha \delta_a) (\phi) =
 \sum_\alpha (-1)^\Abs{\alpha} c_\alpha D^\alpha \phi(a),
\eeq
where the sum extends over all values of the multiindex
$\alpha=(\alpha_1,\ldots,\alpha_d)$ with
$\Abs{\alpha} = \alpha_1+\cdots+\alpha_d \le m$ (the finite degree)
and $D^\alpha = \partial_1^{\alpha_1}\cdots\partial_d^{\alpha_d}$
as usual.

We want to show that $m=0$. Assume the contrary and let $\alpha_0$
be a multiindex that appears in \eqref{point-distr}
with a non-zero coefficient and has $\Abs{\alpha_0} = m$.
Pick a smooth function $\phi$
supported in a neighborhood of $0$ which is such that for each
multiindex $\alpha$ with $\Abs{\alpha} \le m$ we have
$D^\alpha \phi (0) = 0$ if $\alpha \neq \alpha_0$ and
$D^{\alpha_0} \phi(0) = 1$. (To construct such a $\phi$,
multiply the polynomial $1/\alpha_0! x^{\alpha_0}$ with a
smooth function supported in a neighborhood of $0$, which
is identically equal to $1$ in a neighborhood of $0$.)

For $t\to\infty$ let $\phi_t(x) = \phi(t(x-a))$.
Equation \eqref{point-distr} then gives that
\beql{one-side}
\ft{\delta_\Lambda}(\phi_t) = t^m (-1)^m c_{\alpha_0}.
\eeq
On the other hand, using
$$
\FT{\phi(t(x-a))}(\xi) = e^{-2\pi i\inner{a}{\xi/t}} t^{-d}
\ft\phi(\xi/t) 
$$
we get
\beql{other-side}
\ft{\delta_\Lambda}(\phi_t) =
 \sum_{\lambda\in\Lambda} e^{-2\pi i\inner{a}{\lambda/t}} t^{-d}
\ft\phi(\lambda/t).
\eeq
Notice that \eqref{other-side} is a bounded quantity as $t\to\infty$
by a proof similar to that of Theorem \ref{th:at-zero},
while \eqref{one-side} increases like $t^m$, a contradiction.

Hence $\ft{\delta_\Lambda} = w \delta_a$ in a neighborhood of $a$.
The proof of Theorem \ref{th:at-zero} again gives that $\Abs{w} \le
\rho$.
\Qed

Using Theorem \ref{th:point-implies-measure} we may drop from Theorem
\ref{th:periodic-set} the assumption that $\ft{\delta_\Lambda}$ has
to be locallly a measure, as this is now implied by the discrete
support which we assume for $\ft{\delta_\Lambda}$.
Summing up we have the following.
\begin{theorem}\label{th:summary}
Suppose that the multiset $\Lambda$ has uniformly bounded density,
that $S = \supp\ft{\delta_\Lambda}$ is discrete, and that
$$
\Abs{S \cap B_R(0)} \le C R^d,
$$
for some positive constant $C$.
Then $\Lambda$ is a finite union of translated $d$-dimensional
lattices.
\end{theorem}

\mysubsection{Application to tilings by polygons}
\label{sec:polygonal-tilings}

In this section we apply Theorem \ref{th:summary} and the
characterization of the zero-sets of the functions
$\ft{\mu_{e,\tau}}$ (Theorem \ref{th:polygon-zero-set})
in order to give very general
sufficient conditions for a polygon $K$ to admit only
quasi-periodic tilings, if it tiles at all.
\begin{theorem}\label{th:only-quasi-periodic}
Let the polygon $K$ have the pairing property and
tile multiply the plane with the multiset $\Lambda$.
Denote the edges of $K$ by
(we follow the notation of \S\ref{sec:general})
$$
e_1, e_1+\tau_1, e_2, e_2+\tau_2, \ldots, e_n, e_n+\tau_n.
$$
Suppose also that
\beql{no-intersection}
\Set{\widetilde{e_1}, \widetilde{\tau_1}} \cap \cdots \cap 
 \Set{\widetilde{e_n}, \widetilde{\tau_n}} = \emptyset,
\eeq
where with $\widetilde{v}$ we denote the orientation of vector $v$.
Then $\Lambda$ is a finite union of translated
$2$-dimensional lattices.
\end{theorem}
{\bf Proof of Theorem \ref{th:only-quasi-periodic}.}
By Theorem \ref{th:polygon-zero-set} and the tiling assumption
we get
$$
\supp\ft{\delta_\Lambda} \subseteq
 Z(\ft{\mu_{e_1,\tau_1}}) \cap \cdots \cap Z(\ft{\mu_{e_n,\tau_n}}).
$$
By Theorem \ref{th:polygon-zero-set}
in the intersection above each of the sets is contained in
a collection of lines in the direction $\widetilde{e_i}$ union a
collection of
lines 
in the direction $\widetilde{\tau_i}$.
Because of assumption \eqref{no-intersection} these sets
have a discrete intersection
as two lines of different orientations intersect at a point.
Furthermore, because of the regular spacing of these pairs of sets
of lines, it follows that the resulting intersection has at most
$C R^2$ points in a large disc of radius $R$.
Theorem \ref{th:summary} now implies that $\Lambda$ is a finite
union of translated $2$-dimensional lattices.
\Qed

The condition \eqref{no-intersection} is particularly easy to
check for convex polygons.
\begin{theorem}\label{th:convex}
Suppose that $K$ is a symmetric convex polygon which is not a
parallelogram.
Then $K$ admits only quasi-periodic multiple tilings.
\end{theorem}
{\bf Proof of Theorem \ref{th:convex}.}
Suppose that \eqref{no-intersection} fails and that the intersection
in \eqref{no-intersection} contains a vector which is, say, parallel
to the $x$-axis.
It follows that each pair of edges $e_i, e_i+\tau_i$ of edges of $K$
either
(a) has both edges parallel to the $x$-axis, or (b) has the line
joining
the two midpoints parallel to the $x$-axis.
As this latter line goes through the origin it is clear that (b)
can only happen for one pair of edges and, since (a) cannot happen
for two consecutive pairs of edges, (a) can hold at most once as well.
This means that $K$ is a parallelogram.
\Qed

\noindent
{\bf Remarks.}\\
1. It is clear that parallelograms admit tilings which
are not quasi-periodic. Take for example the regular tiling by a
square and move
each vertical column of squares arbitrarily up or down.\\
2. Some very interseting classes of polygons are left out
of reach of Theorem \ref{th:only-quasi-periodic}. An important class
consists of all polygons whose edges are parallel to either the $x$-
or the $y$-axis.\\

\mysection{Bibliography}

\end{document}